%% file: Universality.LiWu.PaperFinal.tex
\input aaWu

\topskip=12pt
\font\sept=cmti9

\def\rightheadline{\ifnum\pageno=\chstart{\hfill}
         \else\centerline
{\sept The universality of symmetric power $\ccmL$-functions
and their Rankin-Selberg $\ccmL$-functions}\hfill
\hskip -3,5mm \tenrm\folio\fi}
\def\leftheadline{\ifnum\pageno=\chstart{\hfill}
         \else\tenrm\folio \hskip -3,5mm \hfill\centerline
{\pauthor H.-Z. Li \& J. Wu}\fi}
\headline={\ifnum\pageno=\chstart{\hfill}
\else{\ifodd\pageno\rightheadline\else\leftheadline\fi}\fi}
\footline={\hfill}

\pageno=1
\newcount\chstart
\chstart=\pageno


\noindent{\it J. Math. Soc. Japan, to appear\hfill}

\vglue 5mm

\centerline{\GGtitre The universality of symmetric power 
$\cmmL$-functions}

\bigskip

\centerline{\GGtitre and their Rankin-Selberg $\cmmL$-functions}

\bigskip

\centerline{{\author H.-Z. Li}
\footnote{$^{(1)}$}{
Project supported by NSFC 10471090, 
the Scientific Research Foundation 
for the Returned Overseas Chinese Scholars 
and Shuguang Plan of Shanghai} (Shanghai) 
\& {\author J. Wu} (Nancy)}

\footnote{{}}{2000 {\it Mathematics Subject Classification}: 11F66}

\footnote{{}}{{\it Key words and phrases}:  
Automorphic $L$-function, Universality}

\vskip 10mm

{\leftskip=1cm
\rightskip=1cm
{\ninerm
{\bf Abstract}.
We establish the universality theorem 
for the first four symmetric power $\ccmL$-functions of automorphic forms
and their associated Rankin-Selberg $\ccmL$-functions.
This generalizes some results of Laurin\u cikas \& Matsumoto
and Matsumoto respectively.
\par}}

\vskip 10mm

\noindent{\bf \S\ 1. Introduction}

\medskip

The automorphic $L$-function is a powerful tool 
to study arithmetic, algebraic and geometric objects. 
Many results will follow from the known or conjectured analytic
properties of automorphic $L$-functions. It is therefore important to
explore an $L$-function in various analytic aspects. Here, we are concerned with the
universality property.  Roughly speaking,
a function $f$ has the universality property if every 
non-vanishing analytic function can be approximated uniformly on compact
subsets in the half critical strip $\D({1\over 2})$ 
by translations of this function $f$ , 
where $\D({1\over 2})$ denotes
$$\D(\sigma_0) 
:= \{s\in \C : \sigma_0<\re s<1\}
\leqno(1.1)$$
for any $\sigma_0<1$. According to Linnik-Ibragimov, it was conjectured
that the universality property is intrinsic to all Dirichlet series
which can be analytically continued to left of their
abscissa of absolute convergence.

The universality of the Riemann zeta-function $\zeta(s)$ 
was first discovered by Voronin [23]. 
More precisely he proved the following:
Let $\K$ be a closed disc of radius $r<{1\over 4}$ 
centered at $s={3\over 4}$,
and $\varphi(s)$ a non-vanishing analytic function in the interior of $\K$
and continuous on $\K$.
Then for any $\varepsilon>0$, there is a real number $t$ such that
$$\sup_{s\in \K}
\big|\zeta(s+it) - \varphi(s)\big|<\varepsilon.
\leqno(1.2)$$
In 1981, Bagchi [1] developed a new method 
to deduce the universality property of $\zeta(s)$ and 
obtained a result stronger than (1.2), as follows.
Let $\K$ be a compact subset of 
\hskip 1mm
$\D({1\over 2})$
with connected complement and $\varphi(s)$ a non-vanishing analytic function in the interior of $\K$
and continuous on $\K$.
Then for any $\varepsilon>0$, we have
$$\liminf_{T\to\infty}
{1\over T}{\rm meas}
\Big\{t\in [0, T] : 
\sup_{s\in \K}
\big|\zeta(s+it) - \varphi(s)\big|<\varepsilon
\Big\}>0,
\leqno(1.3)$$
where ${\rm meas}(\cdot)$ is the Lebesgue measure.
This result was generalized by different authors 
to many other $L$-functions 
such as
Dirichlet $L$-functions,
Dedekind $L$-functions,
Hurwitz $L$-functions,
Lerch $L$-functions, etc.
A detailed historical account can be found in [15].

In this paper we are interested 
in the universality of automorphic $L$-functions.
For a positive even integer $k$ such that $k=12$ or $k\ge 16$
\footnote{$^{(2)}$}{For $k\in \{2, 4, 6, 8, 10, 14\}$,
there is no cusp forms of weight $k$ 
for the full modular group ${\rm SL}(2, \Z)$ (see [21])}, 
we denote by ${\rm H}_k^*$ the set of
all Hecke primitive eigencuspforms of weight $k$ 
for the full modular group ${\rm SL}(2, \Z)$. 
The Fourier series expansion of $f\in {\rm H}_k^*$ 
at the cusp $\infty$ is
$$f(z) = \sum_{n=1}^\infty \lambda_f(n) n^{(k-1)/2} e^{2\pi inz}
\qquad(\im z>0),$$
where $\lambda_f(n)$ is the $n$th (normalized) Fourier coefficient of $f$,
verifying
$$\lambda_f(m)\lambda_f(n)
= \sum_{d\mid (m, n)} \lambda_f\bigg({mn\over d^2}\bigg)
\leqno(1.4)$$
for any integers $m\ge 1$ and $n\ge 1$.
In particular it is a multiplicative function of $n$.
According to Deligne, for any prime number $p$ there is $\alpha_f(p)$
such that
$$\lambda_f(p^\nu)
= \alpha_f(p)^{\nu} + \alpha_f(p)^{\nu-2} + \cdots + \alpha_f(p)^{-\nu}
\qquad(\nu\ge 1)
\leqno(1.5)$$
and
$$|\alpha_f(p)| = 1.
\leqno(1.6)$$
In particular $\lambda_f(1)=1$ and $\lambda_f(n)$ is real.

For $m\in \N$, 
the $m$th symmetric power $L$-function attached to $f\in {\rm H}_k^*$ 
and its Rankin-Selberg $L$-function are defined as
$$L(s, {\rm sym}^mf)
:= \prod_p \prod_{0\le j\le m}
\big(1-\alpha_f(p)^{m-2j} p^{-s}\big)^{-1}
\leqno(1.7)$$
and
$$L(s, {\rm sym}^mf\times {\rm sym}^mf)
:= \prod_p \prod_{0\le i, \, j\le m}
\big(1 - \alpha_f(p)^{2(m-i-j)} p^{-s}\big)^{-1}
\leqno(1.8)$$
for $\sigma>1$, respectively. 
The products over primes in (1.7) and (1.8) admit 
Dirichlet series representation 
$$L(s, F)
= \sum_{n=1}^\infty \lambda_F(n) n^{-s}
\leqno(1.9)$$
for $\sigma>1$,
where $F = {\rm sym}^mf$ or ${\rm sym}^mf\times {\rm sym}^mf$,
and $\lambda_F(n)$ is a multiplicative function. 
Following from (1.6),
we have for $n\ge 1$,
$$|\lambda_F(n)|\le \cases{
\d_{m+1}(n)     & if $F = {\rm sym}^mf$,
\cr\noalign{\medskip}
\d_{(m+1)^2}(n) & if $F = {\rm sym}^mf\times {\rm sym}^mf$,
\cr}
\leqno(1.10)$$
where $\d_{z}(n)$ is the $n$th coefficient of 
the Dirichlet series $\zeta(s)^{z}$.
The case $F = {\rm sym}^1f$ in (1.10) 
is commonly known as Deligne's inequality.

According to [3, Section 3.2.1] and [10, Proposition 2.1],
the gamma factors of $L(s, {\rm sym}^mf)$ 
and $L(s, {\rm sym}^mf\times{\rm sym}^mf)$
are, respectively,
$$L_\infty(s, {\rm sym}^mf) := \cases{
\displaystyle
\prod_{\nu=0}^{n} \Gamma_\C\big(s+(\nu+\dm)(k-1)\big)
& if $m=2n+1$
\cr
\displaystyle
\Gamma_\R(s+\delta_{2\nmid n})
\prod_{\nu=1}^{n} \Gamma_\C\big(s+\nu(k-1)\big)
& if $m=2n$
\cr}
\leqno(1.11)$$
and
$$L_\infty(s, {\rm sym}^mf\times {\rm sym}^mf)
= \cases{\displaystyle
\Gamma_\C(s)^{n+1}
\prod_{\nu=1}^{m} \Gamma_\C\big(s+\nu (k-1)\big)^{m-\nu+1}
& if $m=2n+1$
\cr
\displaystyle
\Gamma_\R(s)
\Gamma_\C(s)^{n}
\prod_{\nu=1}^{m} \Gamma_\C\big(s+\nu (k-1)\big)^{m-\nu+1}
& if $m=2n$,
\cr}
\leqno(1.12)$$
where $\Gamma_\R(s) := \pi^{-s/2}\Gamma(s/2)$,
$\Gamma_\C(s) := 2 (2\pi)^{-s}\Gamma(s)$ and
$$\delta_{2\nmid n} = \cases{
1 & if $2\nmid n$,
\cr\noalign{\smallskip}
0 & otherwise.
\cr}$$
For $F = {\rm sym}^mf$ or $F = {\rm sym}^mf\times {\rm sym}^mf$ where
$f\in {\rm H}_k^*$ and $m=1, 2, 3, 4$, it is known that the function 
$\Lambda(s, F):= L_\infty(s, F) L(s, F)$
is entire on $\C$ and satisfies the functional equation
$$\Lambda(s, F)
= \varepsilon_{F} \Lambda(1-s, F)
\leqno(1.13)$$
with $\varepsilon_{F}=\pm 1$
(see [5, 7, 8, 9] for $F = {\rm sym}^mf$
and [10, 20] for $F = {\rm sym}^mf\times {\rm sym}^mf$).

\smallskip

For the universality property of $L(s, F)$, 
we have the following result.

\proclaim Theorem 1.
Let $1\le m\le 4$,
$2\mid k$ such that $k=12$ or $k\ge 16$, 
$f\in {\rm H}_k^*$ and
$F = {\rm sym}^mf$ or $F = {\rm sym}^mf\times {\rm sym}^mf$.
Define
$$\sigma_F := \cases{
1-(m+1)^{-1} & if $F={\rm sym}^mf$,
\cr\noalign{\medskip} 
1-(m+1)^{-2} & if $F={\rm sym}^mf\times {\rm sym}^mf$.
\cr}
\leqno(1.14)$$
Let $\K$ be a compact subset of 
\hskip 1mm 
$\D(\sigma_F)$ 
with connected complement
and $\varphi(s)$ a non-vanishing analytic function in the interior of $K$
and continuous on $\K$.
Then for any $\varepsilon>0$, we have 
$$\liminf_{T\to\infty}
{1\over T}{\rm meas}
\Big\{t\in [0, T] : 
\sup_{s\in \K}
\big|L(s+it, F) - \varphi(s)\big|<\varepsilon
\Big\}>0.$$

\medskip

{\bf Remark.}
(i)
The particular case $F={\rm sym}^1f = f$ of Theorem 1
was first investigated by 
Ka\u c\.enas \& Laurin\u cikas [6]
and established completely by Laurin\u cikas \& Matsumoto [13].
Another particular case 
$F={\rm sym}^1f\times {\rm sym}^1f = f\times f$ 
was considered by Matsumoto [16] recently.

(ii)
Theorem 1 is established only for $1\le m\le 4$
due to the lack of knowledge about the high symmetric powers.

(iii)
The reason why Theorem 1 holds only for $\D(\sigma_F)$
instead of $\D({1\over 2})$
is that the estimate 
$$\int_1^T \big|L(s, F)\big|^2 \d \tau
\ll_{f} T
\qquad(\forall\;T\ge 1)$$
is only achieved for $\sigma>\sigma_F$ (see (5.3) below),
where $\sigma_F$ is defined as in Theorem 1.
It seems interesting to improve this estimate further
so that Theorem 1 can hold for $\D({1\over 2})$.

(iv)
It is possible to generalize (without too much difficulty) Theorem 1 
to the case of the congruence subgroup $\Gamma_0(N)$ 
with square-free $N$, 
as what Laurin\u cikas, Matsumoto \& Steuding [15] did for $L(s, f)$.

\medskip

Like [13] and [16],
we shall use Bagchi's method to prove Theorem 1. 
(Interested readers are referred to [11] for an excellent paradigm 
on Bagchi's method.)
One of their main tools is Rankin's asymptotic formula
$$\sum_{p\le x} |\lambda_f(p)|^2\sim {x\over \log x}$$
for $x\to \infty$ (see [18], theorem 2).
However, such a prime number theorem for the symmetric $m$th power 
$L$-function with $m\ge 2$ is not available.
In Section 2, we shall establish this result based on [20] and [10],
which is clearly of independent interest 
and may have many other applications.

\smallskip

As in [14], 
we can deduce the following as simple consequences of Theorem 1.

\proclaim Corollary 2. 
Let $m$, $k$, $f$, $F$ and $\sigma_F$ be as in Theorem 1.
For $\sigma_F<\sigma<1$ and any positive integer $J$, 
define a mapping $\psi : \R\to \C^J$ by
$$\psi(\tau)
:= \big(L(\sigma+i\tau, F),
L^{\prime}(\sigma+i\tau, F),
\dots,
L^{(J-1)}(\sigma+i\tau, F)\big).
$$
Then $\psi(\R)$ is dense in $\C^J$.

\proclaim Corollary 3. 
Let $m$, $k$, $f$, $F$ and $\sigma_F$ be as in Theorem 1,
and $J$ be a non negative integer.
If the continuous functions
$g_j : \C^J\to \C\,\,(0\le j\le J)$ satisfy
$$
\sum_{j=0}^J s^j
g_j\big(L(s, F), 
L^{\prime}(s, F),
\ldots,
L^{(J-1)}(s, F)\big)
\equiv 0
$$
for all $s\in \C$, then $g_j\equiv 0$ $(0\le j\le J)$.

\medskip

\noindent{\bf Acknowledgement}.
We began working on this paper in April 2004 
during the visit of the second author 
to Shanghai Jiaotong University,
and finished in February 2005 when the first author visited
l'Universit\'e Henri Poincar\'e (Nancy 1).
We are indebted to both institutions for invitations and support.
The second author would like to thank 
the GRD de Th\'eorie des Nombres au CNRS for support.
We would express our sincere gratitude to K. Matsumoto
for his kind help 
in our study of his joint paper with Laurin\u cikas [13].
Finally the authors wish to thank the referee
for pointing out a mistake in the earlier version.

\vskip 5mm

\noindent{\bf \S\ 2. The prime number theorem  
for symmetric power $L$-functions}

\medskip

Let $m\in \N$, 
$2\mid k$ such that $k=12$ or $k\ge 16$ and $f\in {\rm H}_k^*$.
From (1.6), the product (1.8) is absolutely convergent for $\sigma>1$.
Thus we can define $\Lambda_{{\rm sym}^mf\times {\rm sym}^mf}(n)$
by the relation
$$-{L'\over L}(s, {\rm sym}^mf\times {\rm sym}^mf)
= \sum_{n=1}^\infty 
{\Lambda_{{\rm sym}^mf\times {\rm sym}^mf}(n)\over n^s}
\leqno(2.1)$$
for $\sigma>1$.
The aim of this section is to prove the following result.

\proclaim Proposition 2.1.
Let $1\le m\le 4$, 
$2\mid k$ such that $k=12$ or $k\ge 16$ and $f\in {\rm H}_k^*$.
Then for $x\to\infty$, we have
$$\leqalignno{
& \sum_{n\le x} \Lambda_{{\rm sym}^mf\times {\rm sym}^mf}(n)
\sim x,
& (2.2)
\cr
& \sum_{p\le x} \big|\lambda_f(p^m)\big|^2\log p
\sim x,
& (2.3)
\cr
& \sum_{p\le x} \big|\lambda_f(p^m)\big|^2
\sim {x\over \log x}.
& (2.4)
\cr}$$

This proposition will be referred as the 
{\it prime number theorem}
for the coefficients of symmetric power $L$-functions 
associated with newforms.
It plays a key role in our proof of Theorems 1 \& 2,
and is of independent interest.
The case $m=1$ was considered by Rankin [18].
We shall prove this proposition with the non-vanishing property on
$\sigma=1$ in standard way. To this end, we firstly prove two
preliminary lemmas.

\smallskip

Let $m\in \N$, 
$2\mid k$ such that $k=12$ or $k\ge 16$ and $f\in {\rm H}_k^*$.
We define
$$\Psi_{f, m}(s)
:= \prod_p \prod_{0\le \ell<m}
\big\{\big(1 - \alpha_f(p)^{2(m-\ell)} p^{-s}\big)
\big(1 - \alpha_f(p)^{-2(m-\ell)} p^{-s}\big)\big\}^{-(\ell+1)}
\leqno(2.5)$$
for $\sigma>1$ and
$$\Lambda_{f, m}(n) 
= \cases{
2\sum_{j=1}^m (m+1-j)\cos[2j\theta_f(p)\nu] \log p  
& if $n=p^\nu$,
\cr\noalign{\medskip}
0                                        
& otherwise,
\cr}
\leqno(2.6)$$
where $\alpha_f(p)$ is determined by (1.5)--(1.6)
and $\theta_f(p)\in [0, \pi]$ is chosen such that
$\alpha_f(p) = e^{i\theta_f(p)}$.

\proclaim Lemma 2.1.
Let $m\in \N$, $2\mid k$ such that $k=12$ or $k\ge 16$ and $f\in {\rm H}_k^*$.
Then for $\sigma>1$, we have
$$\leqalignno{-{\Psi_{f, m}'\over \Psi_{f, m}}(s)
& = \sum_{n=1}^\infty {\Lambda_{f, m}(n)\over n^s},
& (2.7)
\cr
\log \Psi_{f, m}(s)
& = \sum_{n=2}^\infty {\Lambda_{f, m}(n)\over n^s \log n}.
& (2.8)
\cr}$$

\noindent{\sl Proof}.
By the Deligne inequality,
the Euler product $\Psi_{f, m}(s)$ converges absolutely for $\sigma>1$.
Taking logarithmic derivative on both sides of (2.5), 
we have, for $\sigma>1$,
$$\eqalign{-{\Psi_{f, m}'\over \Psi_{f, m}}(s)
& = \sum_p \sum_{0\le \ell<m} (\ell+1)
\bigg(
{\alpha_f(p)^{2(m-\ell)} p^{-s}\log p
\over
1-\alpha_f(p)^{2(m-\ell)} p^{-s}}
+ {\alpha_f(p)^{-2(m-\ell)} p^{-s}\log p
\over
1-\alpha_f(p)^{-2(m-\ell)} p^{-s}}
\bigg)
\cr
& = \sum_p \sum_{\nu\ge 1} \sum_{0\le \ell<m} (\ell+1)
{[\alpha_f(p)^{2(m-\ell)\nu}+\alpha_f(p)^{-2(m-\ell)\nu}] \log p
\over p^{s\nu}},
\cr}$$
which is equivalent to (2.7).

Integrating (2.7) on the half-line $\{s+t : t\ge 0\}$, we obtain 
 (2.8).
\hfill
$\square$

\smallskip

\proclaim Lemma 2.2.
Let $1\le m\le 4$, 
$2\mid k$ such that $k=12$ or $k\ge 16$ and $f\in {\rm H}_k^*$.
Then for $\sigma\ge 1$ and $s\not=1$, we have 
$$L(s, {\rm sym}^mf\times {\rm sym}^mf)\not=0.$$

\noindent{\sl Proof}.
Noticing that
$$\sum_{\scriptstyle 0\le i, \, j\le m\atop\scriptstyle i+j=\ell} 1
= \cases{
\ell+1       & if $0\le \ell\le m$
\cr\noalign{\smallskip}
2m-\ell+1    & if $m<\ell\le 2m$
\cr}$$
and using (1.8), we can write, for $\sigma>1$,
$$L(s, {\rm sym}^mf\times {\rm sym}^mf)
= \zeta(s)^{m+1} \Psi_{f, m}(s).
\leqno(2.9)$$

As usual we denote by $\Lambda(n)$ von Mangoldt's function.
Following from (2.6), (2.7), (2.9) and the classical relations
$$-{\zeta'\over \zeta}(s)
= \sum_{n=1}^\infty
{\Lambda(n)\over n^s},
\qquad
\log\zeta(s)
= \sum_{n=2}^\infty
{\Lambda(n)\over n^s\log n}
\qquad(\sigma>1),$$
we infer that
$$\Lambda_{{\rm sym}^mf\times {\rm sym}^mf}(n)
= (m+1)\Lambda(n)
+ \Lambda_{f, m}(n)
\qquad(n\ge 1)
\leqno(2.10)$$
and
$$\log L(s, {\rm sym}^mf\times {\rm sym}^mf)
= \sum_{n=2}^\infty
{\Lambda_{{\rm sym}^mf\times {\rm sym}^mf}(n)\over
n^s\log n}
\qquad(\sigma>1).
\leqno(2.11)$$

Next we calculate $\Lambda_{{\rm sym}^mf\times {\rm sym}^mf}(p^\nu)$.
Write $\vartheta_\nu=\theta_f(p)\nu$ for notational convenience, we get
$$\eqalign{\Lambda_{f, m}(p^\nu)(\log p)^{-1}
& = 2 \sum_{1\le \ell\le m} \sum_{1\le j\le m-\ell+1} 
\cos(\ell 2\vartheta_\nu)
\cr
& = 2 \sum_{1\le j\le m} \sum_{1\le \ell\le m-j+1} 
\cos(\ell 2\vartheta_\nu),
\cr}$$
by (2.6). On the other hand, we have
$$\eqalign{
\sum_{1\le \ell\le m-j+1} \cos(\ell 2\vartheta_\nu)
& = \re\bigg(\sum_{1\le \ell\le m-j+1} e^{i\ell 2\vartheta_\nu}\bigg)
\cr
& = \re\bigg({e^{i(m-j+2)2\vartheta_\nu}-e^{i2\vartheta_\nu}
\over e^{i2\vartheta_\nu}-1}\bigg)
\cr
& = \re\bigg(e^{i(m-j+2)\vartheta_\nu}
{e^{i(m-j+1)\vartheta_\nu}-e^{-i(m-j+1)\vartheta_\nu}
\over e^{i\vartheta_\nu}-e^{-i\vartheta_\nu}}\bigg)
\cr
& = {\cos[(m-j+2)\vartheta_\nu] \sin[(m-j+1)\vartheta_\nu]
\over \sin\vartheta_\nu}. 
\cr}$$
Inserting it into the preceding formula and 
applying the identity 
$$2\cos\alpha\sin\beta=\sin(\alpha+\beta)-\sin(\alpha-\beta),$$
it follows that
$$\eqalign{\Lambda_{f, m}(p^\nu)(\log p)^{-1}
& = {1\over \sin\vartheta_\nu} 
\sum_{1\le j\le m} 
\big(\sin\{[2(m-j+1)+1]\vartheta_\nu\}-\sin\vartheta_\nu\big)    
\cr
& = {1\over \sin\vartheta_\nu} 
\sum_{1\le \ell\le m} \sin[(2\ell+1)\vartheta_\nu]-m.    
\cr}$$
Similarly, we have
$$\eqalign{\sum_{1\le \ell\le m} \sin[(2\ell+1)\vartheta_\nu]
& = \im \bigg(\sum_{1\le \ell\le m} e^{i(2\ell+1)\vartheta_\nu}\bigg)   
\cr
& = \im \bigg(e^{i\vartheta_\nu}
{e^{i(m+1)2\vartheta_\nu}-e^{i2\vartheta_\nu}
\over e^{i2\vartheta_\nu}-1}\bigg)    
\cr
& = \im \bigg(e^{i(m+2)\vartheta_\nu}
{e^{im\vartheta_\nu}-e^{-im\vartheta_\nu}
\over e^{i\vartheta_\nu}-e^{-i\vartheta_\nu}}\bigg)    
\cr
& = {\sin[(m+2)\vartheta_\nu]\sin(m\vartheta_\nu)
\over \sin\vartheta_\nu}    
\cr
& = {\sin^2[(m+1)\vartheta_\nu]-\sin^2\vartheta_\nu
\over \sin\vartheta_\nu}.
\cr}$$
Combining this with the previous relation, we deduce that
$$\Lambda_{f, m}(p^\nu)(\log p)^{-1}
=\bigg({\sin[(m+1)\vartheta_\nu]\over \sin\vartheta_\nu}\bigg)^2-m-1, 
$$
which implies, via (2.10), 
$$\Lambda_{{\rm sym}^mf\times {\rm sym}^mf}(p^\nu)
= \bigg({\sin[(m+1)\vartheta_\nu]\over \sin\vartheta_\nu}\bigg)^2 \log p.
\leqno(2.12)$$
In particular, we obtain with (1.5) that
$$\leqalignno{\Lambda_{{\rm sym}^mf\times {\rm sym}^mf}(p)
& = \bigg({\sin[(m+1)\theta_f(p)]\over \sin\theta_f(p)}\bigg)^2\log p
& (2.13)
\cr
&  = \bigg(\sum_{0\le j\le m} e^{i(m-2j)\theta_f(p)}\bigg)^2\log p
\cr
& = \big|\lambda_f(p^m)\big|^2 \log p.
\cr}$$

Now we are ready to prove Lemma 2.2.
Suppose $L(s, {\rm sym}^mf\times {\rm sym}^mf)$
has a zero at $1+i\tau_0$ of order $\ell\ge 1$, 
where $\tau_0\not=0$.
Consider the function 
$$g(s)
:=L(s, {\rm sym}^mf\times {\rm sym}^mf)^3
L(s+i\tau_0, {\rm sym}^mf\times {\rm sym}^mf)^4
L(s+i2\tau_0, {\rm sym}^mf\times {\rm sym}^mf)^2.$$
Since $L(s, {\rm sym}^mf\times {\rm sym}^mf)$ is holomorphic 
except for a simple pole at $s=1$, 
$g(s)$ is holomorphic for $\sigma\ge 1$ and 
the zero at $s=1$ is of order $\ge 4\ell-3\ge 1$.

But from (2.11), we have for $\sigma>1$,
$$\log g(s)
= \sum_{n\ge 2} 
{\Lambda_{{\rm sym}^mf\times {\rm sym}^mf}(n)\over n^s\log n}
\big(3+4n^{-i\tau_0}+2n^{-i2\tau_0}\big).$$
Together with (2.6), (2.10) and (2.12), we deduce, for $\sigma>1$,
$$\eqalign{\log|g(\sigma)|
& = \sum_{n\ge 2} 
{\Lambda_{{\rm sym}^mf\times {\rm sym}^mf}(n)\over n^\sigma\log n}
\big(3+4\cos(\tau_0\log n)+2\cos(2\tau_0\log n)\big)
\cr
& = \sum_{n\ge 2} 
{\Lambda_{{\rm sym}^mf\times {\rm sym}^mf}(n)\over n^\sigma\log n}
\big(1+2\cos(\tau_0\log n)\big)^2
\cr\noalign{\smallskip}
& \ge 0.
\cr}$$
Thus $|g(\sigma)|\ge 1$ for $\sigma>1$, and 
$g$ cannot have a zero of order $4\ell-3$ $(\ge 1)$ at $\sigma=1$.
This contradiction completes our proof.
\hfill
$\square$

\medskip

Next we shall apply Theorem II.7.11 of [22] 
to prove Proposition 2.1.
Define
$$G(s) 
:= - {L'\over L}(s+1, {\rm sym}^mf\times {\rm sym}^mf) {1\over s+1}
- {1\over s}.$$
Since $\Lambda(s, {\rm sym}^mf\times {\rm sym}^mf)$ is holomorphic 
except for simple poles at $s=0, 1$,
the function $G(s)$ is analytically continued 
to a meromorphic function 
on $\C$.
By Lemma 2.2, we have
$$L(1+i\tau, {\rm sym}^mf\times {\rm sym}^mf)\not=0.$$
Thus $G(s)$ is holomorphic in an open set 
containing the half-plane $\sigma\ge 0$. 
In particular we have
$$\big|G(2\sigma+i\tau)-G(\sigma+i\tau)\big|
\le \sigma 
\sup_{0\le \theta\le 1, \; |\tau|\le T}|G'(\theta+i\tau)|$$
for $T>0$, $0\le \sigma\le \dm$ and $|\tau|\le T$.
From this we deduce
$$\int_{-T}^T \big|G(2\sigma+i\tau)-G(\sigma+i\tau)\big| \d\tau 
= o(1)
\qquad(\sigma\to 0+)$$
for each fixed $T>0$. 
Now Theorem II.7.11 of [22] is applied with 
$F = - L'/L$, $a=c=1$ and $w=0$
to yield the asymptotic formula (2.2).

From (2.13), we can write
$$\sum_{n\le x} \Lambda_{{\rm sym}^mf\times {\rm sym}^mf}(n)
= \sum_{p\le x} \big|\lambda_f(p^m)\big|^2 \log p + R,$$
where we have, via (2.6) and (2.10),
$$\eqalign{R
& := \sum_{p^\nu\le x, \; \nu\ge 2} 
\Lambda_{{\rm sym}^mf\times {\rm sym}^mf}(p^\nu)
\cr
& \,\le \sum_{p\le x^{1/2}} \sum_{\nu\le \log x/\log p}(m+1)^2 \log p
\cr
& \,\le (m+1)^2 \sum_{p\le x^{1/2}} \log x
\cr
& \,\ll_m x^{1/2}.
\cr}$$
Thus (2.2) implies (2.3).
Finally (2.4) follows from (2.3) by integration by parts.
\hfill
$\square$

\vskip 5mm

\noindent{\bf \S\ 3. Bagchi's method and proof of Theorem 1}

\medskip

In this section, 
we present Bagchi's method in our case and first
formulate it as three propositions.
At the end of this section, 
we shall apply Propositions 3.2 and 3.3 to prove Theorem~1.
The proof of these three propositions will be given 
in sections 4, 5 and 6, respectively.

\smallskip

Let $1\le m\le 4$, 
$2\mid k$ such that $k=12$ or $k\ge 16$,
$f\in {\rm H}_k^*$, 
$F = {\rm sym}^mf$ or ${\rm sym}^mf\times {\rm sym}^mf$.
Let $\sigma_F$ and $\D(\sigma_F)$ be defined as in (1.14) and (1.1), 
respectively.
Denote by $H_F$ 
the space of analytic functions on $\D(\sigma_F)$
equipped with the topology of uniform convergence on compact subsets
of $\D(\sigma_F)$.

Let $\gamma:=\{s\in \C : |s|=1\}$ be the unit torus and
$$\Omega := \prod_p \gamma_p,$$
where $\gamma_p = \gamma$ for all prime numbers $p$.
With the product topology and componentwise multiplication,
$\Omega$ is a compact abelian topological group.
Hence there is a unique probability Haar measure $\mu_{\rm h}$
on $(\Omega, {\cal B}(\Omega))$
\footnote{$^{(3)}$}{For any space $X$, 
we denote by ${\cal B}(X)$ the class of all Borel subsets of $X$.}
and we have $\mu_{\rm h}=\prod_p \mu_{{\rm h}, p}$,
where $\mu_{{\rm h}, p}$ is the Haar measure on 
$(\gamma_p, {\cal B}(\gamma_p))$ (see [19], Theorem 5.14). 
For every $\omega = \{\omega_p\}\in \Omega$,
we extend it to a completely multiplicative function, by defining
$$\omega_n := \prod_{p^\nu \| n} \omega_p^\nu.$$
In view of (1.10),
we can prove, similar to Lemma 5.1.6 and Theorem 5.1.7 of [11], 
that there is a subset $\widetilde\Omega\subset \Omega$
with $\mu_{\rm h}(\widetilde\Omega)=1$ such that
for any $\widetilde\omega\in \widetilde\Omega$ 
the series
$$\sum_{n\ge 1} \widetilde\omega_n \lambda_F(n) n^{-s}$$
and the product
$$\prod_p \sum_{\nu\ge 0} 
\widetilde\omega_p^\nu \lambda_F(p^\nu)p^{-\nu s}$$
are uniformly convergent on compact subsets 
of the half-plane $\sigma>\dm$,
and the equality
$$L(s, F; \widetilde\omega)
:= \sum_{n\ge 1} \widetilde\omega_n \lambda_F(n) n^{-s}
= \prod_p \sum_{\nu\ge 0} 
\widetilde\omega_p^\nu \lambda_F(p^\nu)p^{-\nu s}
\leqno(3.1)$$
holds.
Clearly for $\sigma>\dm$ and $\widetilde\omega\in \widetilde\Omega$, 
we have
$$L(s, F; \widetilde\omega)
= \prod_p 
\big(1 + \widetilde\omega_p \lambda_F(p)p^{-s} 
+ O(p^{-2\sigma})\big).$$
Therefore for any $\widetilde\omega\in \widetilde\Omega$, the series
$$\leqalignno{
L^\dagger(s, F; \widetilde\omega) 
& := - \sum_p \log\big(1 - \widetilde\omega_p \lambda_F(p) p^{-s}\big)
& (3.2)
\cr
L_{p_0}^\flat(s, F; \widetilde\omega) 
& := - \sum_{p>p_0} 
\widetilde\omega_p \log\big(1 - \lambda_F(p) p^{-s}\big)
& (3.3)
\cr}$$
are uniformly convergent on compact subsets 
of the half-plane $\sigma>\dm$,
where $p_0\ge 3$ is an arbitrarily fixed constant.
Moreover, we introduce two subsets of $H_F$: 
$$
{\cal L}_F^\dagger 
:= \big\{L^\dagger(s, F; \widetilde\omega) :
\widetilde\omega\in \widetilde\Omega\big\}
\quad{\rm and}\quad
{\cal L}_{F, p_0}^\flat 
:= \big\{L_{p_0}^\flat(s, F; \widetilde\omega) :
\widetilde\omega\in \widetilde\Omega\big\}.
\leqno(3.4)$$

The first auxiliary result of Bagchi's method is 
the denseness of ${\cal L}_F^\dagger$,
which is important in the proof of Proposition 3.3 below.

\proclaim Proposition 3.1.
Let $1\le m\le 4$,
$2\mid k$ such that $k=12$ or $k\ge 16$,
$f\in {\rm H}_k^*$ and
$F={\rm sym}^mf$ or ${\rm sym}^mf\times {\rm sym}^mf$.

\vskip -1,3mm

{\sl 
{\rm (i)}
For any fixed $p_0\ge 3$,
the set ${\cal L}_{F, p_0}^\flat$ is dense in $H_F$.

\smallskip

{\rm (ii)}
The set ${\cal L}_F^\dagger$ is dense in $H_F$.}

\medskip

Define three probability measures on 
$\big(H_F, {\cal B}(H_F)\big)$ 
$$\leqalignno{
P_{F, T}(A)
& := {1\over T} {\rm meas}\big\{t\in [0, T] : L(s+it, F)\in A\big\},
& (3.5)
\cr\noalign{\smallskip}
P_F(A)
& := \mu_{\rm h}\big(\{\omega\in \Omega : 
L(s, F; \omega)\in A\}\big),
& (3.6)
\cr\noalign{\smallskip}
Q_F(A)
& := \mu_{\rm h}\big(\{\omega\in \Omega : 
\log L(s, F; \omega)\in A\}\big),
& (3.7)
\cr}$$
for $A\in {\cal B}(H_F)$.
The next limit theorem is one of the keys of Bagchi's method.

\proclaim Proposition 3.2.
Let $1\le m\le 4$,
$2\mid k$ such that $k=12$ or $k\ge 16$,
$f\in {\rm H}_k^*$
and $F = {\rm sym}^mf$ or ${\rm sym}^mf\times {\rm sym}^mf$.
Then the probability measure $P_{F, T}$ converges weakly to $P_F$
as $T\to \infty$.

\smallskip

The third key step of Bagchi's method is 
to determine the support of the probability measure $P_F$ 
on $(H_F, {\cal B}(H_F))$.
By definition,
a point $s\in S$ is said to be 
in the support of a probability measure $P$
on $(S, {\cal B}(S))$ 
iff every open neighborhood of $s$ has strictly positive measure.
The set of all such points is called the support of $P$, 
denoted by $S(P)$.
Clearly $S(P)$ is the smallest closed subset of $S$ 
such that $P\big(S(P)\big)=1$ (see [2], Chapter 1).
The support of a $S$-valued random variable $Y$ on the probability space
$(X, {\cal B}(X), \mu)$
is the support of the probability measure $P_Y$ on $(S, {\cal B}(S))$
where
$P_Y(A) = \mu(Y\in A)$ ($A\in {\cal B}(S)$),
called the distribution of $Y$.

\proclaim Proposition 3.3.
Let $1\le m\le 4$, 
$2\mid k$ such that $k=12$ or $k\ge 16$,
$f\in {\rm H}_k^*$ and 
$F = {\rm sym}^mf$ or ${\rm sym}^mf\times {\rm sym}^mf$.
With the previous notation, we have the following results:

\vskip -1,3mm

{\sl
{\rm (i)}
The support of the probability measure $Q_F$ 
on $(H_F, {\cal B}(H_F))$ is the whole space $H_F$.

\smallskip

{\rm (ii)}
The support of the probability measure $P_F$ on $(H_F, {\cal B}(H_F))$ is}
$$S_0
:= \big\{\varphi(s)\in H_F : \varphi(s)\not=0 \;
\hbox{for any $s\in \D(\sigma_F)$ or $\varphi(s)\equiv 0$}\big\}.$$

\medskip

Now we apply Propositions 3.1, 3.2 and 3.3 to prove Theorem 1.

\smallskip

Let $\K$ be a compact subset of $\D_\infty(\sigma_F)$
with connected complement.
Let $\varphi(s)$ be a non-vanishing continuous functions on $\K$
which is analytic in the interior of $\K$.
By Lemma 11 of [13],
for any $\varepsilon>0$ we can find a polynomial $p(s)$ 
such that $p(s)\not=0$ on $\K$ and
$$\sup_{s\in \K} |\varphi(s) - p(s)|<\textstyle {1\over 4}\varepsilon.
\leqno(3.8)$$
Since $p(s)$ has only finitely many zeros, 
we can find a region $G_1$ such that 
$\K\subset G_1$ and $p(s)\not=0$ on $G_1$.
We choose $\log p(s)$ to be analytic in the interior of $G_1$.
Applying Lemma 11 of [13] to $\log p(s)$ again, 
we find another polynomial $q(s)$
such that
$$\sup_{s\in \K} 
\big|p(s) - e^{q(s)}\big|<\textstyle {1\over 4}\varepsilon.
\leqno(3.9)$$
From (3.8) and (3.9), we deduce, for any $T>0$,
$$\Big\{t\in [0, T] : 
\sup_{s\in \K} \big|L(s+it, F) - e^{q(s)}\big|
<{\varepsilon\over 2}\Big\}
\subset 
\Big\{t\in [0, T] : 
\sup_{s\in \K} \big|L(s+it, F) - \varphi(s)\big|<\varepsilon\Big\}.
\leqno(3.10)$$
On the other hand, the set 
$$G 
:= \Big\{g\in H_F : 
\sup_{s\in \K}\big|g(s)-e^{q(s)}\big|<\dm \varepsilon\Big\}$$
belongs to $G\in {\cal B}(H_F)$ and is open in $H_F$, thus we have
$$\leqalignno{P_{F, T}(G)
& = {1\over T} {\rm meas}\big\{t\in [0, T] : L(s+it, F)\in G\big\}
& (3.11)
\cr
& = {1\over T} 
{\rm meas}\Big\{t\in [0, T] : 
\sup_{s\in \K}\big|L(s+it, F)-e^{q(s)}\big|<\dm \varepsilon\Big\}.
\cr}$$
By Proposition 3.1, the measure $P_{F, T}(G)$ converges weakly 
to $P_F(G)$ as $T\to\infty$.
With (3.10) and (3.11), Theorem 1.1.8 of [11] leads to 
$$\eqalign{
& \liminf_{T\to\infty}
{1\over T} 
{\rm meas}\Big\{t\in [0, T] : 
\sup_{s\in \K}|L(s+it, F)-\varphi(s)|<\varepsilon\Big\}
\cr
& \ge  \liminf_{T\to\infty}
{1\over T} 
{\rm meas}\Big\{t\in [0, T] : 
\sup_{s\in \K}\big|L(s+it, F)-e^{q(s)}\big|<\dm \varepsilon\Big\}
\cr
& \ge P_F(G).
\cr}$$

Obviously $e^{q(s)}\in S_0 = S(P_F)$ and 
$G$ is a neighbourhood of $e^{q(s)}$.
Therefore $P_F(G)>0$.
This completes the proof of Theorem 1.

\vskip 5mm

\noindent{\bf \S\ 4. Proof of Proposition 3.1}

\medskip

In order to prove Proposition 3.1,
we first apply our result in Section 2 
to establish a preliminary lemma,
which is a generalization of the key lemma 
in Laurin\v cikas \& Matsumoto [13].

\proclaim Lemma 4.1.
Let $1\le m\le 4$,
$2\mid k$ such that $k=12$ or $k\ge 16$ and $f\in {\rm H}_k^*$.
For every $\delta\in [0, 1)$, define
$${\cal P}_\delta 
= {\cal P}_\delta({\rm sym}^mf)
:= \big\{p : 
\hbox{$p$ is prime such that $|\lambda_f(p^m)|\ge\delta$}\big\}.$$
Let $\eta>0$ and $c>1+\eta$ be two fixed constants.
For any $a\ge 2$ and $(1+\eta)a<b\le ca$, we have
$$\sum_{\scriptstyle p\in {\cal P}_\delta\atop\scriptstyle a<p\le b} 
{1\over p}
\ge \bigg\{{1-\delta^2\over (m+1)^2-\delta^2}
+ o_{c, \delta, \eta}(1)\bigg\} 
\sum_{a(1+\eta)<p\le b} 
{1\over p},$$
where $o_{c, \delta, \eta}(1)$ is a quantity 
tending towards 0 as $a\to\infty$.

\noindent{\sl Proof}.
Define 
$$\pi_\delta(x) := \#\big\{p\le x : p\in {\cal P}_\delta\big\}.$$
In particular we have ${\cal P}_0={\cal P}$
(the set of all prime numbers)
and $\pi_0(x)=\pi(x):=\#\big([1, x]\cap {\cal P}\big)$.

Clearly it is sufficient to prove that
for any $(1+\eta)a<u\le b$, 
$$\pi_\delta(u)-\pi_\delta(a)
\ge \bigg\{{1-\delta^2\over (m+1)^2-\delta^2}
+ o_{c, \delta, \eta}(1)\bigg\} 
\big(\pi(u)-\pi(a)\big),
\leqno(4.1)$$
since the desired inequality follows from (4.1) 
via a simple integration by parts.

For $a\le u\le b$, the Deligne inequality $|\lambda_f(p^m)|\le m+1$
allows us to write
$$\leqalignno{\sum_{a<p\le u} \big|\lambda_f(p^m)\big|^2
& \le (m+1)^2 
\sum_{a<p\le u, \, p\in {\cal P}_\delta} 1
+ \delta^2 
\sum_{a<p\le u, \, p\notin {\cal P}_\delta} 1 
& (4.2)
\cr
& \le \big[(m+1)^2-\delta^2\big] 
\big[\pi_\delta(u)-\pi_\delta(a)\big]  
+ \delta^2 \big[\pi(u)-\pi(a)\big].
\cr}$$

According to (2.4) of Proposition 2.1, we have
$$\eqalign{\sum_{a<p\le u} \big|\lambda_f(p^m)\big|^2
& = \pi(u)\{1+o(1)\} - \pi(a)\{1+o(1)\}
\cr
& = \pi(u) - \pi(a) + o\big(\pi(u)\big).
\cr}$$
Since $(1+\eta)a<u\le ca$, a simple calculation shows, 
via the prime number theorem, that
$$\eqalign{\pi(u) - \pi(a)
& = {u\over \log u}\{1+o(1)\} - {a\over \log a}\{1+o(1)\}
\cr
& \ge {\eta a\over \log a}\{1+o_{c, \eta}(1)\}
\cr
& \ge {\eta\over 2c}\pi(u).
\cr}$$
Combining these two estimates yields
$$\sum_{a<p\le u} \big|\lambda_f(p^m)\big|^2
= \big[\pi(u) - \pi(a)\big]\{1+o_{c, \eta}(1)\}.
\leqno(4.3)$$
Now the desired inequality (4.1) follows from (4.2) and (4.3).
\hfill
$\square$

\medskip

Now we are ready to prove Proposition 3.1.

\smallskip

{\it Proof of Proposition 3.1.}
Fix a 
${\widetilde\omega}^0
=\{{\widetilde\omega}^0_p\}
\in \widetilde\Omega$,
then the series
$$L_{p_0}^\flat(s, F; {\widetilde\omega}^0) 
:= - \sum_{p>p_0} {\widetilde\omega}^0_p 
\log\big(1 - \lambda_F(p) p^{-s}\big)
\leqno(4.4)$$
converges in $H_F$.
To prove assertion (i), we shall apply Lemma 4 of [13] 
to this series. In fact, it suffices to verify the condition (a) there,
since conditions (b) and (c) are plainly satisfied.

Let $\mu$ be a complex measure on $(\C, {\cal B}(\C))$
with compact support in $\D(\sigma_F)$ such that
$$\sum_p 
\bigg|\int_\C \log\big(1 - \lambda_F(p) p^{-s}\big) \d\mu(s)\bigg|
<\infty.
\leqno(4.5)$$
Since $\sigma_F>\dm$, we see easily
$$\sum_p |\lambda_F(p)| |\rho(\log p)|
<\infty
\leqno(4.6)$$
with
$$\rho(z) := \int_\C e^{-sz} \d\mu(s).$$

We shall prove that (4.6) leads to
$$\rho(z)\equiv 0,
\leqno(4.7)$$
which implies the validity of condition (a) in Lemma 4 of [13] 
since for any non-negative integer $r$, $\int_\C s^{r} \d\mu(s) = 0$ by differentiating (4.7) $r$-times with respect to $z$ 
and taking $z=0$.
Noticing that 
$$L_{p_0}^\flat(s, F; \widetilde\omega)
= L_{p_0}^\flat\big(s, F; 
(\widetilde\omega/{\widetilde\omega}^0){\widetilde\omega}^0\big)
\quad{\rm and}\quad
\widetilde\omega/{\widetilde\omega}^0
= \big\{\widetilde\omega_p/{\widetilde\omega_p}^0\big\}\in \Omega,$$
Lemma 4 of [13] shows that ${\cal L}_{F, p_0}^\flat$ 
is dense in $H_F$.

It remains to prove (4.7). Firstly we write
$$\rho(z) = \int_\C e^{z s} \d\mu_-(s),$$
where the measure $\mu_-$ is defined by $\mu_-(A)=\mu(-A)$
for $A\in {\cal B}(\C)$ with $-A := \{-a : a\in A\}$.
Clearly $\mu_-$ supports in $\{s\in \C : -1<\sigma<-\dm\}$.
Thus $\rho(z)$ verifies all conditions of Lemma 5 of [13].
If $\rho(z)\not\equiv 0$,
then this lemma implies
$$\limsup_{r\to\infty} {\log |\rho(r)|\over r}
>-1.
\leqno(4.8)$$

Next we shall apply Lemma 7 of [13] 
to deduce an opposite inequality.
This follows a contradiction, and hence (4.7) holds true.

Since the support of $\mu$ is compact and 
is contained in $\D(\sigma_F)$, 
we have
$$|\rho(\pm i y)|
\le e^{My} \int_\C |\d\mu(s)|
\qquad(y>0),$$
where $M=M_\mu$ is a positive constant such that
$\mu$ supports in $(\sigma_F, 1)\times [-M, M]$. 
Thus $\rho(z)$ satisfies condition (a) of Lemma 7 of [13] 
with $\alpha=M$.
Fix a positive number $\beta<\pi/M$,
which assures condition (d) of Lemma 7 of [13].

A similar calculation to (2.13) allows us to obtain
$$\lambda_F(p) = \cases{
\lambda_f(p^m)
& if $F={\rm sym}^mf$,
\cr\noalign{\medskip}
\big|\lambda_f(p^m)\big|^2
& if $F={\rm sym}^mf\times {\rm sym}^mf$.
\cr}
\leqno(4.9)$$

Define
$${\cal L} 
:= \big\{\ell\in \N : \textstyle
\exists \; r\in \big((\ell-{1\over 4})\beta, \; 
(\ell+{1\over 4})\beta\big]
\; \hbox{such that} \; |\rho(r)|\le e^{-r}\big\}.$$ 
By using (4.6) and (4.9), 
we can deduce, for any fixed $\delta\in [0, 1)$,
$$\eqalign{\infty
& > \sum_{p} |\lambda_F(p)| |\rho(\log p)|
\ge \delta^2 \sum_{p\in {\cal P}_\delta} |\rho(\log p)|
\cr\noalign{\smallskip}
& \ge \delta^2 \sum_{\ell\notin {\cal L}} 
\sum_{\scriptstyle p\in {\cal P}_\delta\atop\scriptstyle 
(\ell-{1\over 4})\beta<\log p\le (\ell+{1\over 4})\beta} 
|\rho(\log p)|
\cr
& \ge \delta^2 \sum_{\ell\notin {\cal L}} 
\sum_{\scriptstyle p\in {\cal P}_\delta\atop\scriptstyle 
(\ell-{1\over 4})\beta<\log p\le (\ell+{1\over 4})\beta} p^{-1}. 
\cr}$$
Now we apply Lemma 4.1 with
$$\textstyle
a := \exp\{(\ell-{1\over 4})\beta\},
\qquad
b := \exp\{(\ell+{1\over 4})\beta\},
\qquad
c := \exp\{{1\over 2}\beta\}$$
and $\eta>0$ such that $1+\eta<c$.
It follows that
$$\eqalign{\infty
& > \bigg\{{\delta^2 (1-\delta^2)\over (m+1)^2-\delta^2}
+ o_{c, \delta, \eta}(1)\bigg\} 
\sum_{\ell\notin {\cal L}}
\sum_{(1+\eta)a<p\le b} {1\over p} 
\cr
& \ge \bigg\{{\delta^2 (1-\delta^2)\over (m+1)^2-\delta^2}
+ o_{c, \delta, \eta}(1)\bigg\} 
\sum_{\ell\notin {\cal L}}
\bigg\{\bigg({1\over 2}-{\log(1+\eta)\over \beta}\bigg){1\over \ell}
+ O\bigg({1\over \ell^2}\bigg)\bigg\}, 
\cr}$$
which implies
$$\sum_{\ell\notin {\cal L}} {1\over \ell}
<\infty.
\leqno(4.10)$$
If we write ${\cal L} = \{a_1, a_2, \dots\,\}$
with $a_1<a_2<\dots$, 
it is easy to see that 
$$\lim_{n\to\infty} {a_n\over n} = 1.
\leqno(4.11)$$ 
In fact we have
$[x]=|{\cal L}\cap[1, x]|+|(\N\sset {\cal L})\cap[1, x]|$.
But (4.10) implies
$$|\N\sset {\cal L}\cap[1, x]|
\le \sqrt{x} + \sum_{\sqrt{x}<\ell\le x, \, \ell\notin {\cal L}} 1
\le \sqrt{x} + x\sum_{\ell>\sqrt{x}, \, \ell\notin {\cal L}} 1/\ell
=o(x).$$
Thus $|{\cal L}\cap[1, x]|\sim x$, which is equivalent to (4.11).

By the definition of ${\cal L}$,
there exists a sequence $\{r_n\}$ such that
$$\textstyle
(a_n-{1\over 4})\beta<r_n\le (a_n+{1\over 4})\beta
\qquad{\rm and}\qquad
|\rho(r_n)|\le e^{-r_n}.$$
Then
$$\lim_{n\to\infty} {r_n\over n} = \beta
\qquad{\rm and}\qquad
\limsup_{n\to\infty} {\log |\rho(r_n)|\over r_n}
\le -1.$$
This shows that condition (c) of Lemma 7 of [13] is satisfied.

For any integers $m$ and $n$ such that $m>n\ge 1$, we have
$$r_m - r_n
\ge \textstyle (a_m-a_n-{1\over 2})\beta
\ge \textstyle {1\over 2}(a_m-a_n)\beta
\ge \textstyle {1\over 2}(m-n)\beta.$$
Thus the condition (b) of Lemma 7 of [13] is also satisfied.

Now we can apply Lemma 7 of [13] to write
$$\limsup_{r\to\infty} {\log |\rho(r)|\over r}
= \limsup_{n\to\infty} {\log |\rho(r_n)|\over r_n}
\le -1.$$
This contradicts to (4.8), and
the proof of assertion (i) completes.

Next we shall use the result in assertion (i) to prove (ii). 
Let $\K$ be a compact subset of $\D(\sigma_F)$ and $\varphi\in H_F$.
For any $\varepsilon>0$, we take $p_0\ge 3$ such that
$$\sup_{s\in \K} \sum_{p>p_0} \sum_{\nu\ge 2}
{|\lambda_F(p)|^\nu\over \nu p^{\nu\sigma}}
< {\varepsilon\over 2}.
\leqno(4.12)$$
Since 
$$\varphi(s) + \sum_{p\le p_0}\log\big(1 - \lambda_F(p) p^{-s}\big)
\in H_F,$$
assertion (i) shows that there is 
$\widetilde\omega=\{\widetilde\omega_p\}\in \Omega$ such that
$$\sup_{s\in \K}
\bigg|\varphi(s) 
+ \sum_{p\le p_0}\log\big(1 - \lambda_F(p) p^{-s}\big)
- L_{p_0}^\flat(s, F; \widetilde\omega)
\bigg|
< {\varepsilon\over 2}.
\leqno(4.13)$$
Taking
$$\widetilde\omega_p' := \cases{
1                       & if $p\le p_0$
\cr\noalign{\smallskip}
\widetilde\omega_p      & if $p>p_0$
\cr}
\qquad{\rm and}\qquad
\widetilde\omega' = \{\widetilde\omega_p'\},$$
the inequalities (4.12) and (4.13) imply
$$\eqalign{
\sup_{s\in \K}\big|\varphi(s) 
- L^\dagger(s, F; \widetilde\omega')\big|
& \le \sup_{s\in \K}
\bigg|\varphi(s) 
+ \sum_{p\le p_0}\log\big(1 - \lambda_F(p) p^{-s}\big)
- L_{p_0}^\flat(s, F; \widetilde\omega)\bigg|
\cr
& \quad
+ \sup_{s\in \K} \sum_{p>p_0} 
\big|\log\big(1 - \widetilde\omega_p \lambda_F(p) p^{-s}\big)
-\widetilde\omega_p \log\big(1 - \lambda_F(p) p^{-s}\big)\big|
\cr
& < {\varepsilon\over 2}
+ \sup_{s\in \K} \sum_{p>p_0} \sum_{\nu\ge 2}
{|\lambda_F(p)|^\nu\over \nu p^{\nu\sigma}}
\cr
& < \varepsilon.
\cr}$$
This completes the proof.
\hfill
$\square$

\vskip 5mm

\noindent{\bf \S\ 5. Proof of Proposition 3.2}

\medskip

Obviously Proposition 3.2 is a particular case of Theorem 2 of [12].
Thus it suffices to verify all assumptions there,
that is, to show that there is a positive constant $c$ for which 
$$L(s, F)\ll_f |\tau|^c
\qquad(\sigma>\sigma_F, \; |\tau|\ge 1),
\leqno(5.1)$$
and
$$\int_1^T \big|L(s, F)\big|^2 \d \tau
\ll_{f} T
\qquad(\sigma>\sigma_F, \; T\ge 1).
\leqno(5.2)$$

By using (1.11), (1.12) and (1.13),
a standard Phragm\'en-Lindel\"of argument
allows us to obtain the convex bound for $L(s, F)$, 
i.e. (5.1) with $c=(m+1+\delta_{2\nmid m})/4$ if $F={\rm sym}^mf$ 
and $c=(m+1)^2/4$ if $F={\rm sym}^mf\times {\rm sym}^mf$.
A detailed proof can be found in [10].

In order to verify (5.2), 
we can apply theorem 4 of Perelli [17],
where an estimate of this type was established 
for a general class of $L$-functions.
In view of (1.11) and (1.12),
it is easy to see that $L(s, F)$ lies in the class considered in 
Perelli [17] with evident choice of parameters.
Therefore Theorem 4 of Perelli [17] gives
$$\int_1^T \big|L(s, F)\big|^2 \d \tau
\ll_{f, \varepsilon} T^{(1-\sigma)/(1-\sigma_F)+\varepsilon}
\leqno(5.3)$$
uniformly for $\dm\le \sigma<1$ and $T\ge 1$,
which implies (5.2).
This completes the proof.
\hfill
$\square$

\vskip 5mm

\noindent{\bf \S\ 6. Proof of Proposition 3.3}

\medskip

By the definition,
$\{\omega_p\}$ is a sequence of independent random variables 
defined 
on the probability space $(\Omega, {\cal B}(\Omega), \mu_{\rm h})$,
and the support of each $\omega_p$ is the unit circle $\gamma$.
Hence
$$\bigg\{
\log\bigg(\sum_{\nu\ge 0} 
\omega_p^\nu \lambda_F(p^\nu) p^{-\nu s}\bigg)
\bigg\}$$
is a sequence of independent $H_F$-valued random elements,
and the set
$$\bigg\{\varphi\in H_F :
\varphi(s) =
\log\bigg(\sum_{\nu\ge 0} 
a^\nu \lambda_F(p^\nu) p^{-\nu s}\bigg), \; 
a\in \gamma
\bigg\}$$
is the support of $H_F$-valued random element
$\log\big(\sum_{\nu\ge 0} 
\omega_p^\nu \lambda_F(p^\nu) p^{-\nu s}\big)$.
Consequently,
by theorem 1.7.10 of [11] (see also [13], lemma 10), 
the support of the $H_F$-valued random element
$$\log L(s, F; \omega)
= \sum_p 
\log\bigg(\sum_{\nu\ge 0} 
\omega_p^\nu \lambda_F(p^\nu) p^{-\nu s}\bigg)
$$
is the closure of ${\cal L}^\dagger_F$, 
i.e. the whole space $H_F$ by Proposition 3.1(ii).
This proves the first assertion.

Now we consider any element $\varphi=\varphi(s)$ 
of $S_0\sset \{0\}$ 
and its neighbourhood $G$ in $S_0\sset \{0\}$.
Since the map 
$\exp : H_F\to S_0\sset\{0\}$ is onto and continuous, 
we see that $\exp^{-1}(\varphi)\in H_F$ exists, and 
$\exp^{-1}(G)$ is a neighbourhood of $\exp^{-1}(\varphi)$ 
in $H_F$.   
According to (i), $H_F$ is the support of 
$\log L(s, F; \omega)$, so $Q_F(\exp^{-1}(G))>0$, 
where $Q_F$ is the distribution of $\log L(s,F; \omega)$, 
defined by (3.7).   
But 
$$Q_F(\exp^{-1}(G)) = P_F(G),$$ 
where $P_F$ is the distribution of $L(s, F; \omega)$ 
given by (3.6).   
Hence $P_F(G)>0$.   
This implies that any $\varphi\in S_0\sset\{0\}$ is 
an element of the support of $L(s, F; \omega)$.
Thus
$$S_0\sset\{0\}\subset S(P_F).$$
By Lemma 9 of [13], we have $\overline{S_0\sset\{0\}} = S_0$.
Since $S(P_F)$ is closed, we deduce 
$$S_0\subset S(P_F).
\leqno(6.1)$$

Let $\widetilde\Omega\subset \Omega$ be described as in Section 3.
Then for any $\widetilde\omega\in \widetilde\Omega$,
we have
$$L(s, F; \widetilde\omega)
= \cases{\displaystyle
\prod_p \prod_{0\le j\le m}
\big(1 - \widetilde\omega_p\alpha_f(p)^{m-2j} p^{-s}\big)^{-1}
& if $F={\rm sym}^mf$,
\cr\noalign{\smallskip}
\displaystyle
\prod_p \prod_{0\le i, \, j\le m}
\big(1 - \widetilde\omega_p\alpha_f(p)^{2(m-i-j)} p^{-s}\big)^{-1}
& if $F={\rm sym}^mf\times{\rm sym}^mf$.
\cr}
\leqno(6.2)$$
Since every factor on the right-hand side of (6.1) is non-zero,
the function $L(s, F; \widetilde\omega)$ is also non-vanishing.
Thus
$$\big\{L(s, F; \widetilde\omega) : 
\widetilde\omega\in \widetilde\Omega\big\}
\subset S_0\sset \{0\}$$
and
$$P_F(S_0\sset \{0\})
= \mu_{\rm h}\big(\{\omega\in \Omega : 
L(s, F; \omega)\in S_0\sset \{0\}\big)
\ge \mu_{\rm h}(\widetilde\Omega)
= 1
\quad\Rightarrow\quad
P_F(S_0) = 1.$$
Since $S(P_F)$ is the smallest closed subset of $H_F$ such that 
$P_F\big(S(P_F)\big) = 1$ and $S_0$ is closed, we must have
$$S(P_F)\subset S_0.
\leqno(6.3)$$
Now the required result follows from (6.1) and (6.3).
\hfill
$\square$

\vskip 5mm

\noindent{\bf \S\ 7. Proofs of Corollaries 2 and 3}

\medskip

The proofs of Corollaries 2 and 3 will follow closely those
of Theorems 2 and 3 of [14],
but we reproduce here the details for the convenience of readers.

Let $s_0, \dots, s_{J-1}$ be complex numbers such that $s_0\neq 0$.
Inductively on $J$, we easily see that 
there is a polynomial
$p(s) = \sum_{j=0}^{J-1}b_j s^j$ 
such that
$$
\big(e^{p(s)}\big)^{(j)}\big|_{s=0}
= s_j
\qquad
(0\le j \le J-1).
$$

Let $\sigma_F<\sigma_1<1$, and $\K$ be a compact subset of
$\D(\sigma_F)$ with connected complement such that
$\sigma_1$ is contained in the interior of $\K$. 
We denote by $\delta$ 
the distance of $\sigma_1$ from the boundary of $\K$. 
Then for any $\varepsilon>0$,
Theorem 1 assures that we find a real $\tau$ for which
$$
\sup_{s\in\K} \big|L(s+i\tau,F)-e^{p(s-\sigma_1)}\big|
< {\varepsilon \delta^J\over 2^J J!}
$$
holds. 
Then, using Cauchy's integral formula we have
$$\eqalign{
\big|L^{(j)}(\sigma_1+i\tau,F) - s_j\big|
& = {j!\over 2\pi}
\left|\int_{|s-\sigma_1|=\delta/2}
{L(s+i\tau, F) - e^{p(s-\sigma_1)}
\over (s-\sigma_1)^{j+1}} \d s\right|
\cr
& <\varepsilon
\cr}$$
for $0\le j\le J-1$, which implies Corollary 2.
\hfill
$\square$

\medskip

Next we prove Corollary 3.
Without loss of generality, 
we suppose $g_J\not\equiv 0$. 
Then there exists a bounded region $\G\subset \C^J$ 
and a constant $B_0>0$ such that
$|g_J|\ge B_0$ in $\G$.
 
Let $\sigma\in \D(\sigma_F)$. 
According to Corollary 2, 
we can find a sequence of real
numbers ${\tau_n}\to \infty$ such that 
$$
X_{n}
= \big(L(\sigma+i\tau_n, F),
L'(\sigma+i\tau_n, F),
\dots,
L^{(J-1)}(\sigma+i\tau_n, F)\big)
\in \G.
$$
By the assumption of Corollary 3, 
we have
$$
\sum_{j=0}^{J-1} s^j 
g_j\big(L(s, F), L'(s,F),\ldots,L^{(J-1)}(s,F)\big)
= - s^J 
g_J\big(L(s, F), L^{\prime}(s, F), \dots, L^{(J-1)}(s, F)\big)
$$
for all $s\in \C$. 
Letting $s=\sigma+i\tau_n$ and dividing 
both sides by $(\sigma+i\tau_n)^J$, we obtain
$$
\sum_{j=0}^{J-1} (\sigma+i\tau_n)^{j-J} g_j(X_n)
= - g_J(X_n).
$$
Since $\G$ is bounded, $|g_j(X_n)|$ is bounded ($0\le j\le J-1$). 
Hence the left-hand side of above tends to
zero as $n\to \infty$. 
On the other hand, $|g_J(X_n)|\ge B_0>0$. 
This contradiction finishes the proof of Corollary 3.
\hfill
$\square$

\vskip 8mm

\centerline{\bf References}

\bigskip

\item{[1]}
{\author B. Bagchi}, 
The statistical behaviour and 
universality properties of the Riemann zeta-function 
and other allied Dirichlet series, 
Ph.D. thesis, Indian Statistical Institute, Calcutta, 1981.

\item{[2]}
{\author K.-L. Chung},
{\it A course in probability theory}, 
Third edition, 
Academic Press, Inc., 
San Diego, CA, 2001, xviii+419 pp.

\item{[3]}
{\author J. Cogdell \& P. Michel},
On the complex moments of symmetric power $L$-functions at $s=1$,
{\it IMRN} {\bf 31} (2004), 1561-1618.

\item{[4]}
{\author S. Gelbart \& H. Jacquet},
A relation between automorphic representations of 
$GL(2)$ and $GL(3)$,
{\it Ann. Sci. \'Ecole Norm. Sup. (4)} {\bf 11} (1978), 471-552.

\item{[5]}
{\author H. Iwaniec},
{\it Topics in Classical Automorphic Forms},
Graduate Studies in Mathematics, vol. 17, American Mathematical
Society, Providence, Rhode Island, 1997.

\item{[6]}
{\author A. Ka\v c\v enas \& A. Laurin\v cikas}, 
On Dirichlet series related to certain cusp forms, 
{\it Liet. Mat. Rink.} {\bf 38} (1998), 82--97 (in Russian); 
= {\it Lithuanian Math. J.} {\bf 38} (1998), 64--76. 
MR1663828 (99j:11049)

\item{[7]}
{\author H. Kim},
Functoriality for the exterior square of $GL_4$ 
and symmetric fourth of $GL_2$,
Appendix 1 by Dinakar Ramakrishnan,
Appendix 2 by Henry H. Kim and Peter Sarnak,
{\it J. Amer. Math. Soc.} {\bf 16} (2003), 139--183.

\item{[8]}
{\author H. Kim \& F. Shahidi},
Functorial products for $GL_2\times GL_3$ and
functorial symmetric cube for $GL_2$
(with an appendix by C.J. Bushnell and G. Henniart),
{\it Ann. of Math.} {\bf 155} (2002), 837--893.

\item{[9]}
{\author H. Kim \& F. Shahidi},
Cuspidality of symmetric power with applications,
{\it Duke Math. J.} {\bf 112} (2002), 177--197.

\item{[10]}
{\author Y.-K. Lau \& J. Wu},
A density theorem on automorphic $L$-functions and 
some applications,
{\it Trans. Amer. Math. Soc.} {\bf 358} (2006), 441-472. 

\item{[11]}
{\author A. Laurin\v cikas}, 
{\it Limit Theorems for the Riemann Zeta-Function}, 
Kluwer, Dordrecht, 1996. 

\item{[12]}
{\author A. Laurin\v cikas}, 
On limit distribution of the Matsumoto zeta-function II, 
{\it Liet. Mat. Rink.} {\bf 36} (1996), 464--485 (in Russian); 
= {\it Lithuanian Math. J.} {\bf 36} (1996), 371--387.

\item{[13]}
{\author A. Laurin\v cikas \& K. Matsumoto}, 
The universality of zeta-function attached to certain cusp forms, 
{\it Acta Arith.} {\bf 98} (2001), 345--359.

\item{[14]}
{\author A. Laurin\v cikas \& K. Matsumoto}, 
The joint universality of twisted automorphic L-functions,
{\it J. Math. Soc. Japan}, {\bf 56} (2004), 923--939.

\item{[15]}
{\author A. Laurin\v cikas, K. Matsumoto \& J. Steuding}, 
The universality of $L$-function associated with new forms 
(in Russian),  
{\it Izv. Ross. Akad. Nauk Ser. Mat.} {\bf  67} (2003), no. 1, 83--98;  
translation in 
{\it Izv. Math.} {\bf  67}  (2003), no. 1, 77--90.

\item{[16]}
{\author K. Matsumoto},
The mean values and the universality of Rankin-Selberg $L$-functions,
in: {\it Number theory} (Turku, 1999), 201--221, 
de Gruyter, Berlin, 2001. 

\item{[17]}
{\author A. Perelli},
General $L$-functions,
{\it Ann. Mat. Pura Appl. (4)} {\bf 130} (1982), 287--306.

\item{[18]}
{\author R.A. Rankin},
An $\Omega $-result for the coefficients of cusp forms,
{\it Math. Ann.} {\bf 283} (1973), 239--250.

\item{[19]}
{\author W. Rudin},
{\it Functional Analysis},
Second edition, 
International Series in Pure and Applied Mathematics,
 McGraw-Hill, Inc., New York, 1991, xviii+424 pp.

\item{[20]}
{\author Z. Rudnick \& P. Sarnak},
Zeros of principal $L$-functions and random matrix theory,
{\it Duke Math. J.} {\bf 81} (1996), 269--322.

\item{[21]}
{\author J.-P. Serre}, 
{\it Cours d'arithm\'etiques},
(French) Deuxi\`eme \'edition revue et coorig\'ee.
Le Math\'ematicien, No. 2.
Presses Universitaires de France, Paris, 1977. 188 pp.

\goodbreak

\item{[22]}
{\author G. Tenenbaum}, 
{\it Introduction to analytic and probabilistic number theory},
 Translated from the second French edition (1995) by C. B. Thomas, 
Cambridge Studies in Advanced Mathematics {\bf 46},
Cambridge University Press, 
Cambridge, 1995. xvi+448 pp.

\item{[23]}
{\author S. M. Voronin}, 
A theorem on the "universality" of the Riemann zeta-function, 
{\it Izv. Akad. Nauk SSSR Ser. Mat} {\bf. 39} (1975), 475--486 
(in Russian); 
= {\it Math. USSR-Izv.} {\bf 9} (1975), 443--453.

\vskip 5mm

\hskip -1mm{\author 
Department of Mathematics,
Shanghai Jiaotong University,
1954 Hua Shan Road,
Shanghai 200030,
P. R. of China}

\hskip -1mm{\it E-mail}: {\tt lihz@sjtu.edu.cn}

\medskip

\hskip -1mm{\author Institut Elie Cartan,
UMR 7502 UHP CNRS INRIA,
Universit\'e Henri Poincar\'e (Nancy 1),
54506 Vand\oe uvre--l\`es--Nancy, France}

\hskip -1mm{\it E-mail}: {\tt wujie@iecn.u-nancy.fr}

\end

$$\hskip -90mm
\eqalign{
& \hbox{\author Hongze Li}
\cr\noalign{\vskip -1,5mm}
& \hbox{Department of Mathematics}
\cr\noalign{\vskip -1,5mm}
& \hbox{Shanghai Jiaotong University}
\cr\noalign{\vskip -1,5mm}
& \hbox{1954 Hua Shan Road}
\cr\noalign{\vskip -1,5mm}
& \hbox{Shanghai 200030}
\cr\noalign{\vskip -1,5mm}
& \hbox{P. R. of China}
\cr\noalign{\vskip -1,5mm}
& \hbox{e-mail : \tt lihz@sjtu.edu.cn}
\cr\noalign{\vskip 1,5mm}
& \hbox{\author Jie Wu}
\cr\noalign{\vskip -1,5mm}
& \hbox{Institut Elie Cartan}
\cr\noalign{\vskip -1,5mm}
& \hbox{UMR 7502 UHP-CNRS-INRIA}
\cr\noalign{\vskip -1,5mm}
& \hbox{Universit\'e Henri Poincar\'e (Nancy 1)}
\cr\noalign{\vskip -1,5mm}
& \hbox{54506 Vand\oe uvre--l\`es--Nancy}
\cr\noalign{\vskip -1,5mm}
& \hbox{France}
\cr\noalign{\vskip -1,5mm}
& \hbox{e--mail: \tt wujie@iecn.u-nancy.fr}
\cr}$$

\bye

%% file: aaWu.tex
\input amssym.def
\input amssym.tex

\headline={\ifnum \pageno=1 {\hfill} 
\else{\hss \tenrm -- \folio\ -- \hss}\fi}
\footline={\hfil}

\def\dater{\vglue-10mm\rightline{(\the\day/\the\month/\the\year)}}

\hsize 146mm
\vsize 224mm
\hoffset=6mm
\voffset=8mm
\baselineskip=5mm
\overfullrule =0pt

\font\ninerm=cmr9
 at 10,5pt

\font\GGtitre=cmbx10 at 16pt


\def\og{\leavevmode\raise.30ex
\hbox{$\scriptscriptstyle\langle\!\langle\>$}}    
\def\fg{\leavevmode\raise.24ex
\hbox{$\scriptscriptstyle\>\rangle\!\rangle$}}    

\catcode`\@=11

\font\author=cmcsc10
\font\pauthor=cmcsc10 at 8pt
\font\tenmsx=msam10
\font\sevenmsx=msam10 scaled 700
\font\fivemsx=msam10 scaled 500
\font\tenmsy=msbm10
\font\sevenmsy=msbm10 scaled 700
\font\fivemsy=msbm10 scaled 500
\newfam\msxfam
\newfam\msyfam
\textfont\msxfam=\tenmsx  \scriptfont\msxfam=\sevenmsx
\scriptscriptfont\msxfam=\fivemsx
\textfont\msyfam=\tenmsy  \scriptfont\msyfam=\sevenmsy
\scriptscriptfont\msyfam=\fivemsy

\def\hexnumber@#1{\ifnum#1<10 \number#1\else
\ifnum#1=10 A\else\ifnum#1=11 B\else\ifnum#1=12 C\else
\ifnum#1=13 D\else\ifnum#1=14 E\else\ifnum#1=15 F\fi\fi\fi\fi\fi\fi\fi}

\def\msx@{\hexnumber@\msxfam}
\def\msy@{\hexnumber@\msyfam}
\mathchardef\nmid="3\msy@2D
\mathchardef\varnothing="0\msy@3F
\mathchardef\nexists="0\msy@40
\mathchardef\smallsetminus="2\msy@72
\def\Bbb{\ifmmode\let\next\Bbb@\else
\def\next{\errmessage{Use \string\Bbb\space only in math mode}}\fi\next}
\def\Bbb@#1{{\Bbb@@{#1}}}
\def\Bbb@@#1{\fam\msyfam#1}

\font\tentbl=cmr10 scaled 900
\font\seventbl=cmr7 scaled 900
\font\fivetbl=cmr5 scaled 900

\newfam\tblfam

\textfont\tblfam=\tentbl
\scriptfont\tblfam=\seventbl
\scriptscriptfont\tblfam=\fivetbl


\def \C {{\Bbb C}}
\def \D {{\Bbb D}}

\def \G {{\Bbb G}}

\def \K {{\Bbb K}}
\def \N {{\Bbb N}}

\def \R {{\Bbb R}}

\def \Z {{\Bbb Z}}

\font\ccm=cmmi10 at 9pt
\def\ccmL{\hbox{\ccm L}}

 at 9,5pt

\font\cmm=cmmi10 at 15pt
\def\cmmL{\hbox{\cmm L}}

\def\cmmL{\hbox{\cmm L}}

\font\f=cmr10 at 7pt

\def\f1{\hbox{\f 1}}
\def\f2{\hbox{\f 2}}
\def\f3{\hbox{\f 3}}
\def\f4{\hbox{\f 4}}
\def\f5{\hbox{\f 5}}
\def\f6{\hbox{\f 6}}
\def\f7{\hbox{\f 7}}
\def\f8{\hbox{\f 8}}
\def\f9{\hbox{\f 9}}

\def \d {\,{\rm d}}
\def\re{{\Re e\,}}
\def\im{{\Im m\,}}
\def \dm {{\hbox {${1\over 2}$}}}
\def \sset {{\smallsetminus }}